\documentclass[10pt]{article}
\usepackage{amssymb}

\newcommand{\qed}{$\;\;\;\Box$}
\newenvironment{proof}{\par\smallbreak{\sl Proof.~}}
{\unskip\nobreak\hfill \qed \par\medbreak}

\newcommand{\Con}{C}

\newcommand{\N}{{\mathbb N}}
\newcommand{\R}{{\mathbb R}}
\newcommand{\C}{{\mathbb C}}
\newcommand{\Z}{{\mathbb Z}}
\newcommand{\LL}{{\cal L}}

\renewcommand{\d}{\partial}

\newtheorem{thm}{Theorem}
\newtheorem{lemma}[thm]{Lemma}
\newtheorem{defn}[thm]{Definition}
\newtheorem{cor}[thm]{Corollary}
\newtheorem{rem}[thm]{Remark}

\newcommand{\al}{\alpha}
\newcommand{\be}{\beta}
\newcommand{\ga}{\gamma}
\newcommand{\de}{\delta}
\newcommand{\eps}{\varepsilon}
\newcommand{\vphi}{\varphi}
\newcommand{\la}{\lambda}
\newcommand{\om}{\omega}

\newcommand{\ess}{\mathop{\rm ess}}
\newcommand{\colim}{\mathop{\rm colim}}
\renewcommand{\Re}{\mathop{\mathrm{Re}}\nolimits}

\title{
Linear Hyperbolic Problems
in the Whole Scale \\
of Sobolev-Type Spaces
of Periodic Functions}

\newcounter{thesame}
\author{
I.~Kmit 
\\
{\small
Institute for Applied Problems of Mechanics and Mathematics,}\\
{\small Ukrainian Academy of Sciences}\\
{\small Naukova St.\ 3b,
79060 Lviv,
Ukraine}
\\
{\small   E-mail:
{\tt kmit@informatik.hu-berlin.de}}
}

\date{}

\begin{document}

\maketitle

\begin{abstract}
We study one-dimensional linear hyperbolic systems
with $L^{\infty}$-coeffici\-ents subjected to periodic conditions in time
and reflection boundary conditions in space. 
We derive a priori estimates and give an operator representation of solutions
in the whole scale of Sobolev-type spaces of periodic functions. These spaces
give an optimal regularity trade-off for our problem.
\end{abstract}

\emph{Key words:} hyperbolic systems, periodic-Dirichlet problems,
anisotropic Sobolev spaces, a priori estimates

\emph{MSC 2000: 35L50}

\section{Introduction}\label{sec:intr}

The traveling wave models in 
laser dynamics~\cite{brs,r1,r2,reh,sbw,srs,tro} describe
the dynamical behavior of distributed feedback multisection semiconductor lasers.
The models include
a couple of dissipative semilinear first-order 
hyperbolic equations of a single space variable describing the forward and 
backward propagating complex amplitudes of the light. 
We investigate a linearized version of this
system in the case of small periodic forcing of stationary states (see the discussion in 
Section~\ref{sec:nonlinear}).
Specifically, in the domain $\{(x,t)\,|\,0<x<1$, $-\infty<t<\infty\}$ we
consider system
\begin{equation}\label{eq:1}
\begin{array}{cc}
\partial_tu  + \partial_xu + a(x)u + b(x)v
  = f(x,t)\\
\partial_tv  -  \partial_xv + c(x)u + d(x)v
  = g(x,t),
\end{array} 
\end{equation}
subjected to periodic conditions
\begin{equation}\label{eq:2}
\begin{array}{cc}
u(x,t+T) = u(x,t)\\
v(x,t+T) = v(x,t) \\
\end{array} 
\end{equation}
and reflection boundary conditions
\begin{equation}\label{eq:3}
\begin{array}{cc}
u(0,t) = r_0v(0,t)\nonumber\\
v(1,t) = r_1u(1,t).
\end{array} 
\end{equation}
The unknown functions $u$ and $v$ 
and all the data in~(\ref{eq:1})  are complex 
functions,  $r_0$ and $r_1$ are complex constants, and $T>0$. 
The functions $f$ and $g$ are assumed to be $T$-periodic in $t$ and $a,b,c,d\in L^\infty(0,1)$.

We are motivated by the fact that
a linearization is a first step in many  techniques for local investigation of
nonlinear equations.
Here the linearization is done in a neighborhood of a stationary solution (the coefficients 
in~(\ref{eq:1}) depend only on $x$).  
Our analysis covers practically interesting cases of discontinuous coefficients
and right-hand sides in~(\ref{eq:1}). Discontinuity here correspond to the fact that
different sections of multisection semiconductor lasers have
different electrical and optical properties. 

Our goal is to investigate
uniqueness and solvability questions in  spaces providing us with an optimal
 regularity trade-off (see a preamble of Section~2 for details).
Since~(\ref{eq:1})--(\ref{eq:3}) allows separation of variables, it is natural 
to study our problem within the spaces of periodic in $t$ functions using
their $t$-Fourier representations (see, e.g., \cite[Chapter 5.10]{robinson}
and \cite[Chapter 2.4]{vejvoda}). In Section~2 we introduce
Sobolev-type spaces of periodic functions of any real index $\gamma$ 
including
distributions of any desired order of singularity. They
serve as 
the spaces of solutions to~(\ref{eq:1})--(\ref{eq:3})  and the spaces of right 
hand sides of~(\ref{eq:1}).
We show that, for all sufficiently large $\ga$, the spaces of solutions
are embedded into the algebra $L^{\infty}$ with pointwise multiplication.
This fact makes possible to apply our main results to the aforementioned
nonlinear problems of laser dynamics. In Section~3 we derive 
a priori estimates, thereby proving a uniqueness result. This result is proved under
some not too restrictive and quite applicable to problems
of laser dynamics conditions imposed on the coefficients of~(\ref{eq:1}).
In  Section~4, under some smallness assumption on the coefficients $b$ and $c$, we give an 
operator representation of solutions.

For the clarity of presentation, we restrict ourselves to the $2\times 2$ 
hyperbolic system, which is well enough for applications. 
However,  similar results hold true for the $n\times n$ 
hyperbolic systems (Remark~\ref{rem:nxn}).

In~\cite{KmRe} we consider the case of $b,c\in BV(0,1)$ and $a,d\in L^\infty(0,1)$ and
 prove Fredholm Alternative 
in the Sobolev-type spaces of periodic functions 
continuously embedded into~$L^\infty$.

\section{Sobolev-type spaces of periodic functions and their properties}\label{sec:spaces}

We here construct two scales of Banach spaces $V^\gamma$ (for solutions)
and $W^\gamma$ (for right-hand sides in~(\ref{eq:1})) with a scale 
parameter $\gamma\in\R$, consisting of complex valued functions. 
We will achieve the following properties:
\begin{itemize}
\item
elements of $V^\gamma$ satisfy~(\ref{eq:2}) 
and~(\ref{eq:3}) and elements of $W^\gamma$ satisfy~(\ref{eq:2}), 
\item
elements of
$W^\gamma$ allow discontinuities in $x$, 
\item
 for any $\gamma\in\R$,  the
pair $(V^\gamma,W^\gamma)$ gives an optimal regularity for~(\ref{eq:1})--(\ref{eq:3}). This
means that, from one side, for all $(u,v)\in V^\gamma$ the left-hand side of~(\ref{eq:1}) 
belongs to $W^\gamma$ and, from the other side, for all $(f,g)\in W^\gamma$ solutions
to~(\ref{eq:1})--(\ref{eq:3}) belong to $V^\gamma$.
\end{itemize}
We first introduce Sobolev spaces of $T$-periodic functions 
(see, e.g.~\cite{Garrett,Gorbachuk,robinson,vejvoda}). 
Let $S^T=\R/T\Z$. Define the (Banach) space $\Con^k(S^T)$ of
$k$-times continuously differentiable functions on $S^T$ by
$$
\Con^k(S^T)=\{f: S^T\to\C\,|\,f\circ q\in\Con^k(\R)\},
$$
where $q$ is the quotient map $q: \R\to\R/T\Z$. For the (Freshet) space of smooth functions
we hence have
$$
\Con^{\infty}(S^T)=\bigcap_k\Con^k(S^T).
$$
As a topological vector space this
is a projective limit   of
$\Con^k(S^T)$ with projections being the natural inclusions. Now the space of distributions on $S^T$ is defined as
the ascending union (colimit) of duals of the spaces $\Con^k(S^T)$:
$$
\Con^{\infty}(S^T)^*=\bigcup_k\Con^k(S^T)^*
=\colim\limits_{k\to\infty}\Con^k(S^T)^*.
$$
We are now prepared to  define  Sobolev spaces of periodic functions:
Set $\om=2\pi/T$ and $\vphi_k(t)=e^{ik\om t}$ and define
$$
H^{\ga}(S^T)=\Bigg\{u\in\Con^{\infty}(S^T)^*\,\bigg|\,
\|u\|_{H^{\ga}(S^T)}^2=T^{-1}\sum\limits_{k\in\Z}(1+k^2)^{\gamma}
\left|\left[ u,\vphi_{-k}\right]_{\Con^{\infty}(S^T)}\right|^2<\infty\Bigg\},
$$
where $\left[\cdot,\cdot\right]_{\Con^{\infty}(S^T)}: 
\Con^{\infty}(S^T)^*\times\Con^{\infty}(S^T)\to\C$
is the dual pairing. 

Given $l\in\N_0$, denote
\begin{eqnarray*}
&\displaystyle H^{l,\ga}=H^l(0,1;H^{\ga}(S^T))=\Bigg\{u(\cdot,t): (0,1)\to H^{\ga}(S^T)\,\bigg|\,&
\\&\displaystyle
\|u\|_{H^{l,\ga}}^2=T^{-1}\sum\limits_{k\in\Z}(1+k^2)^{\gamma}
\sum\limits_{m=1}^l\int\limits_0^1
\left|\frac{d^m}{dx^m}\left[ u(x,\cdot),\vphi_{-k}\right]_{\Con^{\infty}(S^T)}
\right|^2\,dx<\infty\Bigg\}.&
\end{eqnarray*}
Set
\begin{equation}\label{eq:Fc}
u_k(x)=T^{-1}\left[ u(x,\cdot),\vphi_{-k}\right]_{\Con^{\infty}(S^T)},
\end{equation}
the $t$-Fourier coefficients of $u\in H^{l,\ga}$.

Finally, for each $\gamma\in\R$ we define the spaces $W^{\gamma}$ and 
$V^{\gamma}$ by
$$
W^{\gamma}=H^{0,\ga}\times H^{0,\ga}
$$
and
$$
V^{\gamma} = \Bigl\{(u,v)\in W^{\gamma}\,|\,
(\d_tu+\d_xu,\d_tv-\d_xv)\in W^{\gamma}\Bigr\}.
$$
These spaces will be endowed with norms
$$
\|(u,v)\|_{W^{\gamma}}^2=\|u\|_{H^{0,\ga}}^2+\|v\|_{H^{0,\ga}}^2
$$
and 
$$
\|(u,v)\|_{V^{\gamma}}^2=\|(u,v)\|_{W^{\gamma}}^2
+\|(\partial_tu+\partial_xu,\partial_tv-\partial_xv)\|_{W^{\gamma}}^2.
$$

Now we describe some useful properties of the function spaces introduced above. 

\begin{lemma}\label{lem:W}
 $W^\gamma$ is a Hilbert space.
\end{lemma}
\begin{proof}
It is known  that $H^\gamma(S^T)$ is a Hilbert space (see, e.g.,~\cite{Garrett}). This
implies that $H^{0,\ga}$ is a Hilbert space as well~\cite{Bourbaki}. The lemma follows. 
\end{proof}

\begin{lemma}\label{lem:complete}
$V^{\gamma}$ is a Banach space.
\end{lemma}

\begin{proof}
Let $(u^{j},v^{j})_{j\in\N}$ be a fundamental sequence 
in $V^{\gamma}$. 
Then $(u^j,v^j)_{j\in\N}$  and $(\d_tu^j+\d_xu^j,\d_tv^j-\d_xv^j)_{j\in\N}$
are fundamental sequences in $W^{\gamma}$. Since $W^{\gamma}$ is complete (Lemma~\ref{lem:W}),
there exist
$(u,v)\in W^{\gamma}$ and $(\tilde u,\tilde v)\in W^{\gamma}$ such that
$$
(u^{j},v^{j})\to (u,v) \; \mbox{ and } \;
(\d_tu^j+\d_xu^j,\d_tv^j-\d_xv^j)\to (\tilde u,\tilde v)
$$
in $W^{\ga}$ as $j\to\infty$.
It remains to show that
$
\d_tu+\d_xu=\tilde u
$
and
$
\d_tv-\d_xv=\tilde v
$
in the sense of generalized derivatives. Indeed, take a smooth 
function 
$\psi: (0,1)\times\left(0,T\right)\to\R$ with compact 
support.
Then
\begin{eqnarray*}
\displaystyle
\left[u,(\d_t+\d_x)\psi\right]=
\lim\limits_{j\to\infty}\left[u^j,(\d_t+\d_x)\psi\right]
=-\lim\limits_{j\to\infty}\left[(\d_t+\d_x)u^j,\psi\right]=-
\left[\tilde u,\psi\right],
\end{eqnarray*}
and similarly with $v$ and $\tilde v$. Here $[\cdot,\cdot]$ is the dual pairing on 
$\Con_0^\infty((0,1)\times (0,T))$. The lemma follows.
\end{proof}

Define
a Euclidian space
$$
E^{\ga}=\Bigg\{(u_k(x))_{k\in\Z}\,\bigg|\,u_k(x)\in L^2(0,1) \mbox{ for\,\,each } k,
\sum\limits_{k\in\Z}
(1+k^2)^{\gamma}\|u_k\|_{L^2(0,1)}^2<\infty\Bigg\}
$$
with inner product
$$
\langle (u_k)_k,(w_k)_k\rangle=\sum\limits_{k\in\Z}(1+k^2)^\gamma\int\limits_{0}^1u_k(x)\overline{w_k(x)}\,dx.
$$

\begin{lemma}
 $E^{\ga}$ is a Hilbert space.
\end{lemma}
\begin{proof}
Let $(u^{j})_{j\in\N}$, where $u^j=(u_k^{j})_{k\in\Z}$ be a fundamental sequence 
in $E^{\ga}$. This means that for any $\eps>0$ there is $N\in\N$ such that
\begin{equation}\label{eq:h1}
\sum\limits_{k\in\Z}(1+k^2)^{\gamma}\int\limits_0^1|u_k^n-u_k^m|^2\,dx<\eps
\end{equation}
for all $n,m\ge N$. It follows that for all $k\in\Z$ the sequence $(u_k^j)_{j\in\N}$ is 
fundamental and hence convergent in $L^2(0,1)$ (by completeness of $L^2(0,1)$). Set
$u_k(x)=\lim_{j\to\infty}u_k^j(x)$ and $u=(u_k)_{k\in\Z}$. Our aim is to show that 
$\sum_{k\in\Z}(1+k^2)^{\gamma}\int\limits_0^1|u_k|^2\,dx<\infty$ and
$\lim_{j\to\infty}u^j(x)=u(x)$ in $E^{\ga}$. Indeed, from~(\ref{eq:h1})
we have
$$
\sum\limits_{|k|\le M}(1+k^2)^{\gamma}\int\limits_0^1|u_k^n-u_k^m|^2\,dx<\eps,
$$
the estimate being uniform in $M\in\N$. Fix $n$ and pass the latter sum to the limit
as $m\to\infty$. We get
$$
\sum\limits_{|k|\le M}(1+k^2)^{\gamma}\int\limits_0^1|u_k^n-u_k|^2\,dx\le\eps,
$$
which is true for any $M\in\N$. This implies that $\lim_{j\to\infty}u^j(x)=u(x)$.
Moreover,
$$
\sum\limits_{k\in\Z}(1+k^2)^{\gamma}\int\limits_0^1|u_k^n-u_k|^2\,dx\le\eps,
$$
Since the series 
$
\sum\limits_{k\in\Z}(1+k^2)^{\gamma}\int_0^1|u_k^n|^2\,dx
$
is convergent for any $n\in\N$, the last inequality implies the convergence of 
the series $\sum_{k\in\Z}(1+k^2)^{\gamma}\int_0^1|u_k|^2\,dx$.
The proof is therewith complete.
\end{proof}

The following lemma is an analog of a result for the Sobolev spaces $H^{\gamma}(S^T)$
in~\cite[Section~2, \S 6]{Gorbachuk}.

\begin{lemma}
The map $u\to (u_k(x))_{k\in\Z}$  is a Hilbert space
isomorphism from
$H^{0,\ga}$ onto~$E^{\ga}$.
\end{lemma}
\begin{proof}
By the definition of $H^{0,\ga}$, for any  $u\in H^{0,\ga}$ the sequence $(u_k(x))_{k\in\Z}$ 
defined by~(\ref{eq:Fc})
 is in $E^{\ga}$. Hence the injectivity of the map $H^{0,\ga}\to E^{\ga}$ 
is a straightforward consequence.

If $\ga\ge 0$, the surjectivity  is a simple,  well-known fact.
Let us prove the surjectivity if $\gamma<0$. We can consider $H^{0,\ga}$ as a space of distributions,
 namely,
$$
 H^{0,\ga}=\left\{u\in L^2(0,1;\Con^{\infty}(S^T))^*\,\Big|\,
\|u\|_{H^{0,\ga}}<\infty\right\}.
$$
Given $(u_k)_{k\in\Z}$ in $E^\ga$, let us define a distribution $u\in  H^{0,\ga}$ by 
$$
[u,f]_{L^2(0,1;\Con^{\infty}(S^T))}=\sum\limits_{k\in\Z}\int\limits_0^1u_{-k}(x)f_{k}(x)\,dx,
$$
where $f\in L^2(0,1;\Con^{\infty}(S^T))$ and $f_{k}(x)$
is the $k$-th Fourier coefficient of $f$ in $t$. The following estimate is straightforward:
\begin{eqnarray*}
\lefteqn{\left|
\sum\limits_{k\in\Z}\int\limits_0^1u_{-k}(x)f_{k}(x)\,dx\right|}
\nonumber\\
&&\le\sum\limits_{k\in\Z}\int\limits_0^1(1+k^2)^{\gamma/2}|u_{-k}(x)|(1+k^2)^{-\gamma/2}|f_{k}(x)|\,dx
\nonumber\\
&&\le
\left(\sum\limits_{k\in\Z}\int\limits_0^1(1+k^2)^{\gamma}|u_k(x)|^2\,dx\right)^{1/2}
\left(\sum\limits_{k\in\Z}\int\limits_0^1(1+k^2)^{-\gamma}|f_{k}(x)|^2\,dx\right)^{1/2}
\nonumber\\
&&=\|(u_k)_{k\in\Z}\|_{E^{\ga}}\|f\|_{H^{0,-\ga}}.
\end{eqnarray*}
This means that $u$ is a continuous linear functional on $H^{0,-\ga}$, namely, that
$u\in L^2(0,1;H^{-\ga}(S^T))^*$.
Since $\Con^{\infty}(S^T)$
is  continuously embedded into $H^{-\ga}(S^T)$, 
$u\in L^2(0,1;\Con^{\infty}(S^T))^*$. The surjectivity of the $t$-Fourier coefficient map $H^{0,\ga}\to E^{\ga}$ 
is therewith proved. Thus, this map is a bijection and it is obviously an isomorphism.
\end{proof}

\begin{cor}\label{cor:id}
For any  $u,v\in H^{l,\ga}$ there exist sequences
 $(u_k)_{k\in\Z},(v_k)_{k\in\Z}$ in $H^l(0,1)$ given by~(\ref{eq:Fc})
such that  the series
\begin{equation}\label{eq:F}
\sum\limits_{k\in\Z}u_k\vphi_k,\qquad
\sum\limits_{k\in\Z}v_k\vphi_k
\end{equation}
converge, respectively, to $u$ and $v$ in $H^{l,\ga}$. Vice versa,
for any sequences  $(u_k)_{k\in\Z}$, $(v_k)_{k\in\Z}$ in $H^l(0,1)$
such that
$\sum_{k\in\Z}(1+k^2)^\gamma\|u_k\|_{H^l(0,1)}^2$ $<\infty$ and 
$\sum_{k\in\Z}(1+k^2)^\gamma\|v_k\|_{H^l(0,1)}^2$ $<\infty$ 
there exist unique $u,v\in H^{l,\ga}$ with $u_k$ and $v_k$
being their $t$-Fourier coefficients.
\end{cor}

In what follows, we will identify distributions $u\in H^{l,\ga}$ and  sequences
$(u_k(x))_{k\in\Z}$ in $H^l(0,1)$ corresponding to these distributions by Corollary~\ref{cor:id}.

\begin{lemma}\label{lem:embed}
$V^{\gamma}$ is continuously embedded into $[H^{1,\ga-1}]^2$.
\end{lemma}

\begin{proof}
Let $(u,v)\in V^{\ga}$. Since $(u,v)\in[H^{0,\ga}]^2$, 
we have $(\d_tu,\d_tv)\in[H^{0,\ga-1}]^2$.
By the definition of~$V^{\ga}$, $(\d_xu,\d_xv)\in[H^{0,\ga-1}]^2$ and
hence $(u,v)\in[H^{1,\ga-1}]^2$. Moreover, we have
\begin{eqnarray*}
&\displaystyle
\|(u,v)\|_{[H^{1,\ga-1}]^2}^2=\|u\|_{H^{0,\ga-1}}^2+\|v\|_{H^{0,\ga-1}}^2+
\|\d_xu\|_{H^{0,\ga-1}}^2+\|\d_xv\|_{H^{0,\ga-1}}^2&\\
&\displaystyle\le
\|u\|_{H^{0,\ga-1}}^2+\|v\|_{H^{0,\ga-1}}^2+\|\d_tu+\d_xu\|_{H^{0,\ga-1}}^2+
\|\d_tv-\d_xv\|_{H^{0,\ga-1}}^2&\\
&\displaystyle
+\|\d_tu\|_{H^{0,\ga-1}}^2+\|\d_tv\|_{H^{0,\ga-1}}^2\le
C\|(u,v)\|_{V^{\gamma}},
\end{eqnarray*}
where the constant $C$ does not depend on $(u,v)$.
\end{proof}

Note that Lemmas~\ref{lem:complete} and~\ref{lem:embed} for $\ga\ge 1$ are proved 
in~\cite{KmRe}. 

\begin{cor}\label{cor:con}
If $\ga>3/2$, then  $V^{\gamma}$ is continuously embedded into $\Bigl[\Con([0,1]\times\R)\Bigr]^2$.
\end{cor}

\begin{proof}
By Lemma~\ref{lem:embed}, $V^{\gamma}\hookrightarrow H^{1}(0,1;H^{\ga-1}(S^T))$ 
continuously. The corollary follows from 
the embedding (see, e.g.,~\cite{herrmann})
$$
H^{\gamma}(S^T)\hookrightarrow\Con(\R),\quad \ga>\frac{1}{2}.
$$
\end{proof}

\begin{cor}\label{cor:trace}
Let $(u,v)\in V^{\ga}$. Then for any $x\in[0,1]$ the traces $u(x,\cdot)$ and
$v(x,\cdot)$ are distributions in $H^{\ga-1}(S^T)$ and
satisfy the estimate
$$
\|(u(x,\cdot),v(x,\cdot))\|_{[H^{\ga-1}(S^T)]^2}^2\le C\|(u,v)\|_{V^{\ga}}^2,
$$
where $C$ does not depend on $x$, $u$, and $v$.
\end{cor}
\begin{proof}
 The corollary follows from  
the continuous embedding
$$
V^{\gamma}\hookrightarrow H^{1}(0,1;H^{\ga-1}(S^T))\hookrightarrow \Con(0,1;H^{\ga-1}(S^T)).
$$
\end{proof}

\section{A priori estimates}\label{sec:uniq}

We here give  conditions ensuring the uniqueness of generalized solutions 
to~(\ref{eq:1})--(\ref{eq:3}).
We start from the definition of a generalized solution.

\begin{defn}\label{defn:gen}
A function $(u,v)\in V^{\gamma}$ is called a generalized
solution to the problem~(\ref{eq:1})--(\ref{eq:3}) if it satisfies~(\ref{eq:1})
in $H^{0,\ga}$ and~(\ref{eq:3}) in $H^{\ga-1}(S^T)$.
\end{defn}

To formulate the main result of this section, we will make the following assumption
about the coefficients of the differential equations and the reflection coefficients $r_0$ and $r_1$:
There exist $p,q\in\R$ such that one of the following conditions is fulfilled:
\begin{equation}\label{eq:coef+}
\begin{array}{cc}
\displaystyle
\ess\inf\Bigl[2\Re a-|b|^{2p}-|c|^{2q}\Bigr]>0
\\[2mm]\displaystyle
\ess\inf\Bigl[2\Re d-|b|^{2(1-p)}-|c|^{2(1-q)}\Bigr]>0
\end{array}
\end{equation}
or
\begin{equation}\label{eq:coef-}
\begin{array}{cc}
\displaystyle
\ess\inf\biggl[-2\Re a-\bigg|\frac{b}{r_0}\bigg|^{2p}-\bigg|\frac{c}{r_1}\bigg|^{2q}\biggr]
+2(1-m)\ln\frac{1}{|r_0r_1|}>0
\\[3mm]\displaystyle
\ess\inf\biggl[-2\Re d-\bigg|\frac{b}{r_0}\bigg|^{2(1-p)}-
\bigg|\frac{c}{r_1}\bigg|^{2(1-q)}\biggr]+
2m\ln\frac{1}{|r_0r_1|}>0,
\end{array}
\end{equation}
for $m=0$ or  $m=1$.

\begin{thm}\label{thm:apriori}
Let $\ga\in\R$ be a fixed real, $a,b,c,d\in L^{\infty}(0,1)$, and $(f,g)\in W^{\gamma}$. 
Assume that  $|r_0|<1$ and $|r_1|<1$.
If one of the conditions (\ref{eq:coef+}),~(\ref{eq:coef-}) with $m=0$,
and~(\ref{eq:coef-}) with $m=1$ is true,
then every generalized solution to the problem~(\ref{eq:1})--(\ref{eq:3})
satisfies the  a priory estimate
\begin{equation}\label{eq:+}
\|(u,v)\|_{V^{\gamma}}\le C\|(f,g)\|_{W^{\gamma}}
\end{equation}
for some $C>0$ not depending on $(f,g)$.
\end{thm}

Note that the assumptions $|r_0|<1$ and $|r_1|<1$ are caused by physical reasons.

\begin{proof}
Due to the assumptions imposed on the functions $f$ and $g$, they allow the following
series representations:
\begin{equation}\label{eq:fg0}
\sum\limits_{k\in\Z}f_k\vphi_k,\qquad
\sum\limits_{k\in\Z}g_k\vphi_k,
\end{equation}
where $f_k(x)=T^{-1}\left[ f(x,\cdot),\vphi_{-k}\right]_{\Con^{\infty}(S^T)}$ and 
$g_k(x)=T^{-1}\left[ g(x,\cdot),\varphi_{-k}\right]_{\Con^{\infty}(S^T)}$. Clearly,
$f_k, g_k\in L^2(0,1)$,
\begin{equation}\label{eq:fg}
\sum\limits_{k\in\Z}(1+k^2)^{\gamma}\|f_k(x)\|_{L^2(0,1)}^2<\infty,\quad
\sum\limits_{k\in\Z}(1+k^2)^{\gamma}\|g_k(x)\|_{L^2(0,1)}^2<\infty,
\end{equation}
and the series~(\ref{eq:fg0}) converge to $f$ and $g$ in $H^{0,\gamma}$.
Assume that $(u,v)$ is a generalized solution to the
problem~(\ref{eq:1})--(\ref{eq:3}).
Represent $u$ and $v$ as series~(\ref{eq:F}). Hence
$u_k,v_k$ for each $k\in\Z$ are in $H^1(0,1)$ and satisfy the boundary value problem 
\begin{equation}\label{eq:ordinary}
\begin{array}{cc}
\displaystyle
u_k^{\prime}=f_k(x)-(a(x)+ik\omega)u_k-b(x)v_k\\
v_k^{\prime}=-g_k(x)+(d(x)+ik\omega)v_k+c(x)u_k,
\end{array}
\end{equation}
\begin{equation}\label{eq:two}
\begin{array}{cc}
\displaystyle
u_k(0)=r_0v_k(0)\\
v_k(1)=r_1u_k(1).
\end{array}
\end{equation}
Our aim is to show that 
\begin{equation}\label{eq:7.1}
\sum\limits_{k\in\Z}(1+k^2)^{\gamma}\Bigl[\|u_k(x)\|_{L^2(0,1)}^2+
\|v_k(x)\|_{L^2(0,1)}^2\Bigr]<\infty,
\end{equation}
\begin{equation}\label{eq:u_kv_k}
\sum\limits_{k\in\Z}(1+k^2)^{\gamma}\Bigl[\|u_k^{\prime}(x)+ik\omega u_k(x)\|_{L^2(0,1)}^2+
\|v_k^{\prime}(x)-ik\omega v_k(x)\|_{L^2(0,1)}^2\Bigr]<\infty.
\end{equation}
The estimate~(\ref{eq:u_kv_k}) follows from~(\ref{eq:fg}),~(\ref{eq:7.1}), and~(\ref{eq:ordinary}).
It remains to prove (\ref{eq:7.1}).  We will distinguish two cases.

{\it Case 1. Condition~(\ref{eq:coef+}) is fulfilled.}
Fix $k\in\Z$. Multiplying the
equations of the system~(\ref{eq:ordinary}) by $\overline u_k$ and $\overline v_k$,
respectively, and then summing up the resulting equalities
with their complex conjugations, we  arrive at
the system
\begin{equation}\label{eq:UN+}
\begin{array}{cc}
\displaystyle
\int\limits_0^1\Bigl(\overline u_ku_k^{\prime}+u_k\overline u_k^{\prime}\Bigr)\,dx+
2\int\limits_0^1\Re a|u_k|^2\,dx
\\\displaystyle
=\int\limits_0^1\Bigl(f_k\overline u_k+\overline f_ku_k\Bigr)\,dx
-\int\limits_0^1\Bigl(b\overline u_kv_k+\overline bu_k\overline v_k\Bigr)\,dx
\\\displaystyle
\int\limits_0^1\Bigl(\overline v_kv_k^{\prime}+v_k\overline v_k^{\prime}\Bigr)\,dx-
2\int\limits_0^1\Re d|v_k|^2\,dx
\\\displaystyle
=-\int\limits_0^1\Bigl(g_k\overline v_k+\overline g_kv_k\Bigr)\,dx
+\int\limits_0^1\Bigl(c\overline v_ku_k+\overline cv_k\overline u_k\Bigr)\,dx.
\end{array}
\end{equation}
Using integration by parts and  boundary conditions~(\ref{eq:two}), we get
$$
\int\limits_0^1\Bigl(\overline u_ku_k^{\prime}+u_k\overline u_k^{\prime}\Bigr)\,dx=
|u_k(1)|^2-|r_0|^2|v_k(0)|^2,
$$
$$
\int\limits_0^1\Bigl(\overline v_kv_k^{\prime}+v_k\overline v_k^{\prime}\Bigr)\,dx=
|r_1|^2|u_k(1)|^2-|v_k(0)|^2.
$$
Subtraction of  the second equality of~(\ref{eq:UN+}) from the first one yields
\begin{equation}\label{eq:UN+-}
\begin{array}{cc}
\displaystyle
(1-|r_1|^2)|u_k(1)|^2+(1-|r_0|^2)|v_k(0)|^2
+2\int\limits_0^1\Bigl(\Re a|u_k|^2+
\Re d|v_k|^2\Bigr)\,dx
\\\displaystyle
=\int\limits_0^1\Bigl(f_k\overline u_k+\overline f_ku_k\Bigr)\,dx+
\int\limits_0^1\Bigl(g_k\overline v_k+\overline g_kv_k\Bigr)\,dx
\\\displaystyle
-\int\limits_0^1\Bigl(b\overline u_kv_k+\overline bu_k\overline v_k\Bigr)\,dx
-\int\limits_0^1\Bigl(c\overline v_ku_k+\overline cv_k\overline u_k\Bigr)\,dx.
\end{array}
\end{equation}
Since  $|r_0|<1$ and $|r_1|<1$,
the sum of the  boundary terms is positive.
We will make  use of the following simple inequalities: Given $\eps_0>0$, we have
$$
\biggl|\int\limits_0^1f_k\overline u_k\,dx\biggr|\le
\frac{1}{2\eps_0}\int\limits_0^1|f_k|^2\,dx+
\frac{\eps_0}{2}\int\limits_0^1|u_k|^2\,dx,
$$
$$
\biggl|\int\limits_0^1b\overline u_kv_k\,dx\biggr|\le
\int\limits_0^1|b|^p|u_k||b|^{1-p}|v_k|\,dx,
$$
$$
\le\frac{1}{2}\|b\|_{L^{\infty}(0,1)}^{2p}\int\limits_0^1|u_k|^2\,dx+
\frac{1}{2}\|b\|_{L^{\infty}(0,1)}^{2(1-p)}\int\limits_0^1|v_k|^2\,dx,
$$
 We estimate all other integrals
in the right-hand side of~(\ref{eq:UN+-}) similarly. Finally, by
assumption~(\ref{eq:coef+}), one can choose $\eps_0>0$ so small that
\begin{equation}\label{eq:U0}
\|u_k(x)\|_{L^2(0,1)}^2+
\|v_k(x)\|_{L^2(0,1)}^2
\le C
\Bigl[\|f_k(x)\|_{L^2(0,1)}^2+\|g_k(x)\|_{L^2(0,1)}^2\Bigr],
\end{equation}
where the constant $C$ depends on $a,b,c,d$, but not on $k$. 
Now,~(\ref{eq:fg}) and~(\ref{eq:U0}) imply~(\ref{eq:7.1}).

{\it Case 2. Condition~(\ref{eq:coef-}) is fulfilled.} We first give the proof under
the condition~(\ref{eq:coef-}) with $m=0$. 
We start from the observation that  $(u_k,v_k)$ is a solution
to the problem~(\ref{eq:ordinary})--(\ref{eq:two}) iff
 $(w_k,v_k)$, where $e^{x\al+(1-x)\be}w_k=u_k$ and $\al,\be\in\R$ are fixed reals, is a solution to
the problem
\begin{equation}\label{eq:ordinary1}
\begin{array}{cc}
\displaystyle
w_k^{\prime}=e^{-\al x-(1-x)\be}f_k(x)-(a(x)+ik\omega+\al-\be)w_k
-e^{-\al x-(1-x)\be}b(x)v_k\\[2mm]
v_k^{\prime}=-g_k(x)+(d(x)+ik\omega)v_k+c(x)w_ke^{\al x+(1-x)\be},
\end{array}
\end{equation}
\begin{equation}\label{eq:two1}
\begin{array}{cc}
\displaystyle
e^{\be}w_k(0)=r_0v_k(0)\\[2mm]
v_k(1)=r_1e^{\al}w_k(1).
\end{array}
\end{equation}
 Let us write down an analog of the equality~(\ref{eq:UN+-})
for the problem~(\ref{eq:ordinary1})--(\ref{eq:two1}):
\begin{eqnarray*}
&\displaystyle
-\Bigl(1-|r_1|^2e^{2\al}\Bigr)|w_k(1)|^2-\Bigl(1-|r_0|^2e^{-2\be}\Bigr)|v_k(0)|^2&
\\&\displaystyle
+2\int\limits_0^1\Bigl[(-\Re a-\al+\be)|w_k|^2-
\Re d|v_k|^2\Bigr]\,dx&
\\&\displaystyle
=-\int\limits_0^1e^{-x\al-(1-x)\be}\Bigl(f_k\overline w_k+\overline f_kw_k\Bigr)\,dx-
\int\limits_0^1\Bigl(g_k\overline v_k+\overline g_kv_k\Bigr)\,dx&
\\&\displaystyle
+\int\limits_0^1e^{-x\al-(1-x)\be}\Bigl(b\overline w_kv_k+\overline bw_k\overline v_k\Bigr)\,dx
+\int\limits_0^1e^{x\al+(1-x)\be}\Bigl(c\overline v_kw_k+\overline cv_k\overline w_k\Bigr)\,dx.&
\end{eqnarray*}
Fix $\eps_1>0$ and set $\al=\ln|r_1|^{-1}+\eps_1$ and $\be=-\ln|r_0|^{-1}-\eps_1$.
This clearly forces
$$
-(1-|r_1|^2e^{2\al})|w_k(1)|^2-(1-|r_0|^2e^{-2\be})|v_k(0)|^2>0.
$$
On the  account of the assumption~(\ref{eq:coef-}),
 similarly to  Case~1, $\eps_1>0$ can be chosen so small
that estimate~(\ref{eq:U0})  with $u_k$ replaced by $w_k$ is true
for some $C$  independent of~$k$. Replacing 
$w_k$ with $e^{-x\al-(1-x)\be}u_k$, we arrive at the estimate~(\ref{eq:U0}).

The estimate~(\ref{eq:U0}) under the condition~(\ref{eq:coef-}) with $m=1$ can be obtained
in much the same way, the only difference being in considering the
problem~(\ref{eq:ordinary})--(\ref{eq:two}) with
$v_k=e^{x\al+(1-x)\be}w_k$. 

The estimate~(\ref{eq:7.1}) is proved.
This finishes the  proof of the theorem.
\end{proof}
The following corollary  is straightforward.

\begin{cor}\label{uniq}
Under the conditions of Theorem~\ref{thm:apriori}
a generalized solution to the problem~(\ref{eq:1})--(\ref{eq:3})
(if such exists) is unique.
\end{cor}

\section{Operator representation of generalized solutions}~\label{sec:explicit}

In this section, under some smallness assumption on the coefficients $b$ and $c$, we give an 
explicit formula (an operator representation) 
of a generalized solution. 

Set
$$
V^{\ga}(r_0,r_1)=\Bigl\{(u,v)\in V^{\ga}\,\Big|\,u(0,\cdot)=r_0v(0,\cdot),\ v(1,\cdot)=r_1u(1,\cdot)
\Bigr\},
$$
where the traces $u(0,\cdot),u(1,\cdot),v(0,\cdot)$, and $v(1,\cdot)$  
are interpreted  as distributions in $H^{\ga-1}(S^T)$ according to
Corollary~\ref{cor:trace}.
Given $a, b, c, d\in L^{\infty}$, let us introduce linear operators $A\in\LL(V^{\gamma}(r_0,r_1);W^{\gamma})$ by
$$
A\left[
\begin{array}{c}
u\\v
\end{array}
\right]
=
\left[
\begin{array}{c}
\partial_tu+\partial_xu+au\\\partial_tv-\partial_xv+dv
\end{array}
\right]
$$
and $B\in\LL(W^{\gamma})$ by
$$
B\left[
\begin{array}{c}
u\\v
\end{array}
\right]
=
\left[
\begin{array}{c}
bv\\cu
\end{array}
\right].
$$

\begin{lemma}\label{lem:iso}
Assume that $a, b, c, d\in L^{\infty}(0,1)$ and
\begin{equation}\label{eq:coef}
|r_0r_1|\ne \exp{{\int\limits_0^1(\Re a+\Re d)\,dx}}.
\end{equation}
Then the operator $A$ is
an isomorphism from $V^{\gamma}(r_0,r_1)$ onto $W^{\gamma}$.
\end{lemma}

In the case $\ga\ge 1$  this lemma is proved in~\cite{KmRe}. 

\begin{proof}
Fix an arbitrary $(f,g)\in W^{\ga}$. The functions $f$ and  $g$ are represented by the
series~(\ref{eq:fg0})  with coefficients $f_k(x)\in L^2(0,1)$ and  $g_k(x)\in L^2(0,1)$
satisfying~(\ref{eq:fg}).
We have to show that there exists exactly one $(u,v)\in V^{\ga}(r_0,r_1)$ such that
$$
\d_tu+\d_xu+au=f,\quad \d_tv-\d_xv+du=g.
$$
Representing $u$ and $v$ as the series~(\ref{eq:F}),  we have to show 
that there exists exactly one pair of sequences $(u_k)_{k\in\Z}$ and $(v_k)_{k\in\Z}$
with  $u_k,v_k\in H^1(0,1)$ satisfying~(\ref{eq:two}),~(\ref{eq:7.1}),~(\ref{eq:u_kv_k}), and
\begin{equation}\label{eq:3.3}
\begin{array}{cc}
\displaystyle
u_k^{\prime}+(a(x)+ik\omega)u_k=f_k(x)\\
v_k^{\prime}-(d(x)+ik\omega)v_k=-g_k(x).
\end{array}
\end{equation}

To simplify the formulae below, let us introduce the following notation:
$$
\al(x)=\int\limits_0^xa(y)\,dy,\quad \de(x)=\int\limits_0^xd(y)\,dy,
$$
$$
\Delta_k=e^{ik\om+\de(1)}-r_0r_1e^{-ik\om-\al(1)}.
$$
By a straightforward calculation, the boundary value problem~(\ref{eq:3.3}),~(\ref{eq:two})
has a unique solution $(u_k,v_k)\in [H^1(0,1)]^2$, and this solution is explicitely given by
\begin{equation}\label{eq:3.7}
\begin{array}{rcl}
\displaystyle
u_k(x)&=&e^{-ik\om x-\al(x)}\left(
\displaystyle\int\limits_0^xe^{ik\om y+\al(y)}f_k(y)\,dy
+\frac{r_0}{\Delta_k}w_k(f_k,g_k)\right),\\
v_k(x)&=&e^{ik\om x+\de(x)}\left(
\displaystyle\int\limits_0^xe^{-ik\om y-\de(y)}g_k(y)\,dy
+\frac{1}{\Delta_k}w_k(f_k,g_k)\right)
\end{array}
\end{equation}
with
$$
w_k(f_k,g_k)=r_1e^{-ik\om-\al(1)}
\int\limits_0^1e^{ik\om y+\al(y)}f_k(y)\,dy-
e^{ik\om+\de(1)}
\int\limits_0^1e^{-ik\om y-\de(y)}g_k(y)\,dy.
$$
Here we used assumption~(\ref{eq:coef}), which implies
\begin{equation}\label{eq:3.9}
|\Delta_k|\ge\left|e^{\de(1)}-|r_0r_1|e^{-\al(1)}\right|>0
\end{equation}
for all $k\in\Z$. From~(\ref{eq:3.7}) and~(\ref{eq:3.9})  it follows that
\begin{equation}\label{eq:3.10}
|u_k(x)|+|v_k(x)|\le C\Bigl(\|f_k\|_{L^2(0,1)}+\|g_k\|_{L^2(0,1)}\Bigr)
\end{equation}
for all $x\in[0,1]$, where the constant $C$ does not depend on $k$, $f_k$, $g_k$, and $x$.
Finally,~(\ref{eq:fg}) and~(\ref{eq:3.10}) imply~(\ref{eq:7.1}).
Noting that the estimate~(\ref{eq:u_kv_k}) follows from~(\ref{eq:fg}),~(\ref{eq:7.1}), and~(\ref{eq:3.3}), 
we finish the proof.
\end{proof}

\begin{thm}\label{thm:explicit}
Let $\ga$ be a real,  $a, b, c, d\in L^{\infty}(0,1)$,
 and $(f,g)\in W^{\gamma}$.
Assume that  
the condition~(\ref{eq:coef})
is true. Suppose also that
\begin{equation}\label{eq:norm}
\begin{array}{c}
\displaystyle
\biggl[1+\exp{\left\{3\|a\|_{L^{\infty}}+3\|d\|_{L^{\infty}}\right\}}
\left(1+(1+|r_0|)(1+|r_1|)\left|e^{\de(1)}-|r_0r_1|e^{-\al(1)}\right|^{-1}\right)
\\
\displaystyle
\times
\left(1+\|a\|_{L^{\infty}}+\|d\|_{L^{\infty}}\right)\biggr]
\Bigl(\|b\|_{L^{\infty}}+\|c\|_{L^{\infty}}\Bigr)<1.
\end{array}
\end{equation}
Then the problem~(\ref{eq:1})--(\ref{eq:3})
has a unique generalized solution given by the formula
\begin{equation}\label{eq:explicit}
(u,v)=\sum\limits_{n=0}^{\infty}(-A^{-1}B)^nA^{-1}(f,g).
\end{equation}
\end{thm}

\begin{proof}
By Lemma~\ref{lem:iso}, the problem~(\ref{eq:1})--(\ref{eq:3})  is equivalent to
\begin{equation}\label{eq:01}
(u,v)=-A^{-1}B(u,v)+A^{-1}(f,g),
\end{equation}
where $A^{-1}$ is defined by means of~(\ref{eq:3.7}).
Since  $\|A^{-1}B\|_{\LL(V^{\ga};V^{\ga})}$
 is bounded from above by the left-hand side of~(\ref{eq:norm}), 
$\|A^{-1}B\|_{\LL(V^{\ga};V^{\ga})}<1$. 
Since $V^{\ga}$ is a Banach space for any $\ga\in\R$ (Lemma~\ref{lem:complete}), application of the
 Banach fixed point theorem 
to the equation~(\ref{eq:01}) gives the unique solvability of the latter.
Hence~(\ref{eq:1})--(\ref{eq:3}) has a unique generalized solution.
Iteration of~(\ref{eq:01}) now gives the desired formula~(\ref{eq:explicit}).
\end{proof}

\begin{rem}\label{rem:nxn}
The results of Sections~3 and~4 (Theorems~\ref{thm:apriori} and~\ref{thm:explicit})
can be easily generalized to $n\times n$ hyperbolic systems, namely, to the 
problems of the following kind:
$$
\begin{array}{cc}
\displaystyle
\partial_tu_j  + \lambda_j(x)\partial_xu_j +\sum\limits_{k=1}^na_{jk}(x)u_k = f_j(x,t),
\quad j=1,\dots,n,\quad 0<x<1,\quad t\in\R\\
\displaystyle
u_j(x,t+T) = u_j(x,t),\qquad j=1,\dots,n,\quad 0<x<1,\quad t\in\R\\
\displaystyle
u_j(0,t) = \sum\limits_{k=m+1}^nr_{jk}^0u_k(0,t),\qquad j=1,\dots,m,\quad t\in\R\\
\displaystyle
u_j(1,t) = \sum\limits_{k=1}^mr_{jk}^1u_k(1,t),\qquad j=m+1,\dots,n,\quad t\in\R.
\end{array}
$$
\end{rem}

\section{Concluding remarks}\label{sec:nonlinear}
Let us turn back to the area of laser dynamics that is a motivation of this paper. 
Our overall goal here, which remains a subject of future research, is
to obtain a local existence and uniqueness result for semilinear hyperbolic systems
with small periodic forcing
\begin{equation}\label{eq:nl}
\begin{array}{cc}
\partial_tu  + \partial_xu + g_1(x,u,v,\lambda)
  = \eps f_1(x,t,u,v,\la,\eps)\\
\partial_tv  -  \partial_xv + g_2(x,u,v,\lambda)
  = \eps f_2(x,t,u,v,\la,\eps).
\end{array} 
\end{equation}
Let us choose a solution space to be $V^\gamma$ with $\gamma>3/2$, which is 
advantageous due to Corollary~\ref{cor:con}. A natural way to achieve our goal
can now consist in application of the Implicit Function Theorem to~(\ref{eq:nl}).
For this purpose we would need to have, first, the isomorphism property of a linearization of~(\ref{eq:nl}),
second, the $\Con^1$-smoothness property of the Nemitsky composition operators defined by
the nonlinearities of~(\ref{eq:nl}).

Explicit sufficient conditions for the former property are provided by the results 
of~\cite{KmRe} and Theorem~\ref{thm:apriori} of this paper
(note that the explicitness here is really important from the point of view of
applications). More specifically, assume that $(u,v)=(0,0)$ is a (stationary) solution
to~(\ref{eq:nl}) and linearize this system in a neighborhood of the
stationary solution. Putting $\eps=0$ and $\lambda=0$, we arrive at
the system~(\ref{eq:1}) with zero right hand side. Denote  the operator
corresponding to this system by $F$.
In~\cite{KmRe} we proved that, if the condition~(\ref{eq:coef}) is fulfilled, then $F$ is a Fredholm operator 
from $V^\gamma$ onto $W^\gamma$. 
Theorem~\ref{thm:apriori} states constructive
conditions ensuring the injectivity of $F$. Since any Fredholm operator is an
isomorphism between two Banach spaces iff it is injective, we therefore have the desired isomorphism property for $F$
under rather wide  explicit conditions on the data of~(\ref{eq:nl}).
Note that Theorem~\ref{thm:explicit} ensures the isomorphism property for a 
range of data, in which we cannot use the Fredholmness result from~\cite{KmRe}.

\subsection*{Acknowledgments}
I thank  Lutz Recke for helpful discussions.

\end{document}